\documentclass[12pt]{amsart}
\usepackage{cite,amssymb,graphicx,caption}
\usepackage{amsfonts}
\usepackage{enumerate}
\usepackage{epsfig}
\usepackage{color}
\usepackage{psfrag}
\usepackage{hyperref}
\usepackage{babel}
\usepackage{lineno}
\usepackage{chemarr}
\usepackage{rotating}
\usepackage{graphicx}

\allowdisplaybreaks

\usepackage{rotating}

\newcommand{\R}{\ensuremath{\mathbb{R}}}

\newtheorem{lemma}{Lemma}

\def\eps{\epsilon}
\def\ds{\displaystyle}

\def\rme{\mathrm{e}}
\def\rmi{\mathrm{i}}

\def\criteps{\eps_{\mathrm{c}}}
\def\endeps{\eps_{\mathrm{end}}}

\usepackage{pgfplots}
\usepackage{stackengine}
\usepgfplotslibrary{fillbetween}
\pgfplotsset{compat=newest}
\usetikzlibrary{decorations.markings, patterns}

\usepackage[short,nodayofweek,level,12hr]{datetime}

\def\TL{\mathrm{TL}}
\def\TA{\mathrm{TA}}
\def\CP{\mathrm{CP}}
\def\JS{\mathrm{JS}}

\title[Scaling in stochastic ghosts]{Semiclassical theory explains stochastic ghosts scaling}

\author{J. Tom\'as L\'azaro}
\address{$^{\TL}$Departament de Matem\`atiques (Universitat Polit\`ecnica de Catalunya UPC), Av. Diagonal 647, 08028 Barcelona, Spain}
\address{$^{\TL,\TA,\JS}$Centre de Recerca Matem\`atica (CRM). Edifici C, Campus de Bellaterra, 08193 Cerdanyola del Vall\`es, Barcelona, Spain}
\address{$^{\TL}$Institute of Mathematics of the UPC-BarcelonaTech (IMTech),  Pau Gargallo 14, 08028 Barcelona, Spain}
\address{$^{\TL,\TA,\JS}$Dynamical Systems and Computational Virology, CSIC Associated Unit CRM-Institute for Integrative Systems Biology (I$^2$SysBio), Spain.}
\vspace{3pt}

\author{Tom\'as Alarc\'on}
\address{$^{\TA}$ICREA, Pg. Lluis Companys 23, 08010 Barcelona, Spain}
\address{$^{\TA}$Departament de Matem\`atiques, Universitat Aut\`onoma de Barcelona, Barcelona, Spain}
\vspace{3pt}

\author{Carlos P. Garay}
\address{$^{\CP}$Laboratorio Subterr\'aneo de Canfranc (LSC), 22880, Canfranc-Estaci\'on, Huesca, Spain}
\address{$^{\CP}$Institute for Integrative Systems Biology (I$^2$SysBio), CSIC-UVEG, 46071, Valencia, Spain.}
\vspace{3pt}

\author{Josep Sardany\'es}
\vspace{3pt}

\begin{document}

\maketitle

\begin{abstract}
Slowing down phenomena occur in both deterministic and stochastic dynamical systems at the vicinity of phase transitions or bifurcations. 
An example is found in systems exhibiting a saddle-node bifurcation, which undergo a dramatic time delay towards equilibrium. Specifically the duration of the transient, $\tau$, close to this bifurcation in deterministic systems follows scaling laws of the form $\tau \sim |\epsilon - \epsilon_c|^{-1/2}$, where $\epsilon$ is the bifurcation or control parameter, and  $\epsilon_c$ its critical value. For systems undergoing a saddle-node bifurcation, the mechanism involves transients getting trapped by a so-called ghost. In a recent article we explored how intrinsic noise affected the deterministic picture. Extensive numerical simulations showed that, although scaling behaviour persisted in the presence of noise, the scaling law was more complicated than a simple power law. In order to gain deeper insight into this scaling behaviour, we resort to the WKB asymptotic approximation of the Master Equation. According to this approximation, the behaviour of the system is given as the weighted sum of \emph{trajectories} within the phase space of the Hamiltonian associated to the corresponding Hamilton-Jacobi equation. By analysing the flight time of the Hamilton equations, we show that the statistically significant paths follow a scaling function that exactly matches the one observed in the stochastic simulations. We therefore put forward that the properties of the flight times of the Hamiltonian system underpin the scaling law of the underlying stochastic system, and that the same properties should extend in a universal way to all stochastic systems whose associated Hamiltonian exhibits the same behaviour.
\end{abstract}

\vspace{2pc}
\noindent{\it Keywords}: Bifurcations,  Complex systems, First-order phase transitions, Scaling laws, Stochastic dynamics, Transients, Universality. 
%\submitto{\NJP}
%\maketitle

%\tableofcontents

\section{Introduction}

Slowing down phenomena are ubiquitous in dynamical systems, both stochastic and deterministic, approaching a phase transition or a bifurcation. Both these terms refer to a qualitative change in the dynamical behaviour of the system as a parameter, referred to as control or bifurcation parameter, reaches its critical value \cite{Goldenfeld1992,Kuznetsov1998,Strogatz2000}. In the vicinity of bifurcation points, slowing down phenomena occur whereby the transients towards equilibrium become much longer than those experienced far away from the critical bifurcation point \cite{Hohenberg1977,Strogatz1989,Strogatz2000,Suzuki1982I,Suzuki1982II}. 

An specific example can be found in deterministic systems undergoing a saddle-node (s-n) bifurcation. In this bifurcation, the equilibria collide and get annihilated, as the colliding equilibria become complex after crossing the bifurcation~\cite{Strogatz2000,Sardanyes2006,Fontich2008,Sardanyes2019,Gimeno2018}. Recent research has explored the dynamical mechanisms taking place at the complex phase space causing such slowdown in transients~\cite{Canela2022}. These delays are usually referred to as delayed transitions or ghosts~\cite{Strogatz1989,Strogatz2000}. The term ghost alludes to the attracting remnant left by the s-n collision in the region of the phase space where it occurs. This collision produces a bottleneck which hugely delays the dynamics of the system. Ghosts have been identified in mathematical models of charge density waves~\cite{Strogatz1989}, hypercycles~\cite{Sardanyes2006,Sardanyes2007}, and ecological systems with facilitation including semi-arid ecosystems~\cite{Vidiella2018} and metapopulations with~\cite{Sardanyes2019} and without~\cite{Fontich2008,Gimeno2018} habitat destruction. 

A remarkable property of  slowing down behaviour close to critical points is that it is characterised by universal scaling laws \cite{Hohenberg1977,Strogatz1989,Strogatz2000,Suzuki1982I,Suzuki1982II}. For the s-n bifurcation, the scaling behaviour manifests in the form of a power law of the relaxation time, $\tau$, towards a stable equilibrium. Specifically $\tau \sim |\epsilon - \epsilon_c|^{-1/2}$, where $\epsilon$ is the control parameter and $\epsilon_c$ its critical value~\cite{Strogatz2000,Fontich2008,Duarte2011}. Interestingly, this scaling behaviour has been clearly identified experimentally in an electronic circuit~\cite{Trickey1998}.

Ghosts have been thoroughly analysed in deterministic systems~\cite{Strogatz1989,Trickey1998,Strogatz2000,Sardanyes2006,Sardanyes2007,Fontich2008,Kuehn2008,Sardanyes2008,Sardanyes2010,Vidiella2018,Canela2022}.  Recent computational studies have addressed the effects of intrinsic noise on the ghost phenomenon~\cite{Sardanyes2006,Sardanyes2007,Vidiella2018}. In a recent paper~\cite{SardanyesNJP2020}, we investigated by means of extensive numerical simulations the effects of intrinsic noise on delayed transitions by using a simple model for autocatalysis and a two-member hypercycle. Somewhat contrary to intuition, we found that the ghost phenomenon was robust to intrinsic noise. Indeed, we showed that in stochastic systems close to a first-order phase transition, the relaxation dynamics towards the (unique) absorbing state experienced a slowing down when it entered into the region of the phase space where the collision of (deterministic) equilibria occurred, and therefore the ghost phenomenon was robust to intrinsic noise. We quantified this property using the average extinction time, $\bar{T}_E$. Extensive computational analysis of the $\bar{T}_{E}$ allowed us to ascertain that it exhibited a scaling behaviour.  While the extinction time in deterministic systems close to the s-n bifurcation exhibits simple power law scaling, $\bar{T}_E$ followed a more complex scaling behaviour of the form 
\begin{equation}\label{eq:TE}
\bar{T}_E\sim\Omega^{-b}\,{\mathcal{G}}(\Omega^a(\epsilon-\epsilon_c)),
\end{equation}
where $\Omega$ was the system size, $\epsilon$ and $\epsilon_c$ were the control parameter and its critical value (with $\epsilon\gtrsim\epsilon_c$), respectively, $a$ and $b$ were scaling exponents, and ${\mathcal{G}}(\cdot)$ the scaling function.  The form of the scaling function ${\mathcal{G}}$ is shown in the upper panel of Fig.~\ref{fig_1}.

A theory explaining the shape of the scaling function ${\mathcal{G}}$ is currently lacking. The aim of the present work is to advance this issue by providing a theoretical framework to understand the shape of ${\mathcal{G}}$. Here, we address this problems from the perspective of the Hamiltonian formalism associated with the WKB approximation of the Master Equation (ME) and its path integral representation. By analysing the trajectories of the corresponding Hamiltonian system and their statistical weights of two different models with a s-n bifurcation, as prescribed by the path integral solution of the ME, we show that the shape of ${\mathcal{G}}$ is determined by the properties of the Hamiltonian dynamical system (HDS) close to the bifurcation point. In particular, we study the flight times of the HDS in the vicinity of the bifurcation i.e. as a function of $\phi\equiv\epsilon-\epsilon_c$, and the statistical weight of the corresponding trajectories. This theory allows us to describe in detail the features of the scaling function. Furthermore, since the behaviour of interest depends only on a few features of the Hamiltonian phase space, we claim that this behaviour is universal to all the stochastic systems which share those features. 

The article is organised as follows. In Section 2 we summarise the theoretical framework used to analyse the stochastic dynamics. Sections 3 and 4 provide a detailed account of our results and their discussion, respectively.

\section{Theoretical framework}

In this section we summarise how intrinsic noise is introduced in the general context of nonlinear birth-and-death processes. Assuming that the underlying stochastic dynamics is described by a Markovian process, the basic description of the process is given by the corresponding Master Equation:

\begin{equation} \label{eq:ME}
\frac{\partial}{\partial t} P(X, t) = \sum_{r_i=\pm 1} \left(W_i(X - r_i)P(X-r_i) - W_i(X)P(X)\right),
\end{equation}

\noindent where $X(t)$ is the state variable which corresponds to the number of individuals at time $t$, $P(X,t)$ is the probability density that the system has $X$ individuals at time $t$, and the transition rates, $W_i$, and the stoichiometric coefficients, $r_i$, are such that:

$$\mbox{Prob}\left(X(t+\Delta t)=X+r_i\vert X(t)=X\right)=W_i(X)\Delta t+{\mathcal{O}}(\Delta t^2).$$

In order to proceed further with our analysis, we assume that the transition rates, upon rescaling, satisfy the following scaling relation:

\begin{equation}\label{eq:scalingrates}
W_i(X)=\Omega w_{i}(x)+{\mathcal{O}}(\Omega^0)
\end{equation}

\noindent where $\Omega$ is the so-called system size and $x=X/\Omega$ is the rescaled state variable. Note that rates derived from Law of Mass Action kinetics and Michaelis-Menten-Hill kinetics satisfy such conditions. If the transition rates satisfy Eq. (\ref{eq:scalingrates}), we can propose a Wentzel-Kramers-Brillouin (WKB) Ansatz for $P(X,t)$:

$$P(X,t)=\exp\left(-\Omega(S(x,t)+{\mathcal{O}}(\Omega^{-1}))\right),$$

\noindent which, when introduced in Eq. (\ref{eq:ME}), at the lowest order (i.e. at ${\mathcal{O}}(\Omega^0)$), the function $S(x,t)$, the so-called action, satisfies the equation

\begin{equation}\label{eq:hamiltonjacobi}
\frac{\partial S}{\partial t}=H\left(x,\frac{\partial S}{\partial x}\right)
\end{equation}

\noindent where the function $H(x,p)$ is given by:

\begin{equation}\label{eq:hamiltonian}
H(x,p)=\sum_{r_i=\pm 1}\left(\rme^{r_ip}-1\right)\omega_i(x). 
\end{equation}

The action also maps the stochastic process to a path integral representation\cite{Kubo1973}, formalism that can be studied using tools developed in equilibrium statistical physics. The path integral replaces the deterministic unique trajectory by a functional integral over an infinity of stochastic possible trajectories weighted by the action along each trajectory ($\rme^{-\Omega S}$). The path integral provides a mapping U$_t$ between probability generating functions at different times and a field theory in the continuum limit. The transition probability can be written as 
\begin{equation}
P(x_f,t_f|x_i,t_i)=\iint  dx\,dp\,\rme^{ -\Omega \int(p \dot{x} - H(x,p))\, dt } 
\label{eq:pathint}
\end{equation}
where all possible trajectories in the phase space are weighted by the action, which satisfies a Hamilton-Jacobi equation with Hamiltonian given by Eq.~\eqref{eq:hamiltonian}. This suggests a classical mechanical interpretation, where $H(x,p)$ determines the motion of a \emph{particle} in a phase space where $x$ is the position and $p$ the conjugate momentum. The trajectories of such a particle are given by the solution of the associated Hamilton's equations:
\begin{equation}
\begin{array}{rcrcl} 
\dot{x}&=&{\ds \frac{\partial H}{\partial p}} &=&{\ds \sum r_i \rme^{r_ip}w_i(x)}, \\[1.2ex]
\dot{p}&=&{\ds -\frac{\partial H}{\partial x}} &=&{\ds \sum (\rme^{r_ip}-1)\frac{dw_i}{dx}}. 
\label{eq:hamiltoneqsp} 
\end{array}
\end{equation}
The Hamiltonian formulation of the stochastic processes given by Eqs.~\eqref{eq:hamiltoneqsp} has been used to study such problems as the quasi-steady state distribution in connection with first passage problems \cite{hinch2005,assaf2006,escudero2009,Assaf2010,gottesman2012,bressloff2014a}, the statistics of rare events in reaction-diffusion systems \cite{elgart2004}, and the path integral representation of the process \cite{Kubo1973}, among others. Here, we use the Hamiltonian formulation, specifically the structure of the corresponding phase space, to determine which of its generic features lend robustness to the mean-field ghosts. We further investigate how the scaling form of the relaxation time is an emergent property of such generic features.

\section{Results}

Using the theoretical framework discussed in the previous section, we are ready to explore stochastic systems whose mean-field limit exhibit a s-n bifurcation. Our goal is to explain the numerically observed robustness of the ghost to intrinsic noise and the scaling of the mean extinction times, $\bar{T}_E$, in the autocatalytic replicator model $\dot{x} = k \, x^2 (1 - x) - \epsilon \, x$ (see~\cite{SardanyesNJP2020} for the stochastic analysis and~\cite{Sardanyes2007,Fontich2008,Sardanyes2019} for further information of this model). For completeness we also analyse the Hill model~\cite{Hill1910} with linear decay, which also exhibits the same mean-field bifurcation, and is given by $\dot{x} = k \, x^2/(A^2 + x^2) - \epsilon \, x$, having $k, \epsilon, A > 0$ for both systems. We anticipate that the behaviour of these two models is similar close to the bifurcation and thus our results will be presented alternating both systems. For simplicity, we will hereafter denote Gillespie simulations as GS.
%TABLE %%%%%%%%%%%%%%%%%%%%%%%%%%%%%%%%%%%%%%%%%%%%%%%%%%%%%%%%%%%%%

\begin{table}[]
\centering
\begin{tabular}{|l|l|l|l|}
\hline
{\bf Model } & {\bf Transition}    & {\bf Stoichiometry} & {\bf Description} \\ 
 {\bf system }            &         {\bf rates}                           & {\bf  coefficient} &  \\ 
\hline
Autocatalyic  & $W_1=k_1X$                                 & $r_1=+1$        & Birth       \\ \hline
Autocatalyic  & $W_2=\frac{k_1}{C \Omega^2}X(X-1)(X-2)$    & $r_2=-1$       & Competition for finite resources \\ \hline
Autocatalyic  & $W_3=\epsilon X$                           & $r_3=-1$       & Death       \\ \hline
Hill      & $W_1=\Omega\frac{k_1X^2}{\Omega^2A^2+X^2}$ & $r_1=+1$           & Birth with saturation      \\ \hline
Hill          & $W_2=\epsilon X$                           & $r_2=-1$      & Death       \\ \hline
\end{tabular}
\caption{Transition rates for the two specific model examples used in this paper: the autocatalytic model~\cite{SardanyesNJP2020} and the Hill system~\cite{Hill1910} with exponential decay.} \label{tab:table_1}
\end{table}
%%%%%%%%%%%%%%%%%%%%%%%%%%%%%%%%%%%%%%%%%%%%%%%%%%%%%%%%%%%%%%%%%%%%%%
%FIGURE %%%%%%%%%%%%%%%%%%%%%%%%%%%%%%%%%%%%%%%%%%%%%%%%%%%%%%%%%%%%
\begin{figure*}
\begin{center}
\centerline{\includegraphics[width=0.68\textwidth]{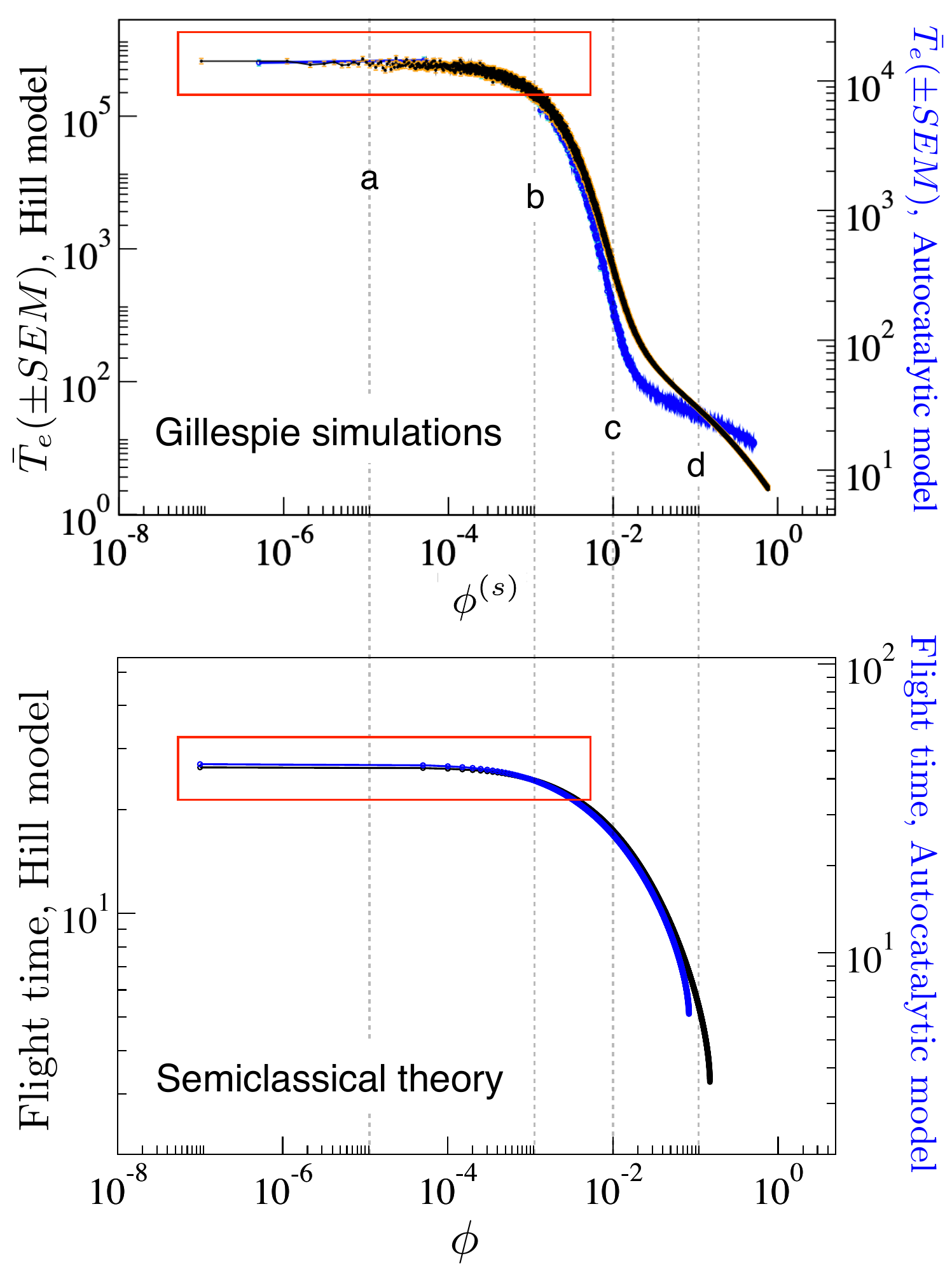}}
\caption{(Top) Scaling of the mean extinction times $\bar{T}_e (\pm SEM)$ (computed averaging over $10^3$ independent realisations) obtained from Gillespie simulations, shown overlapped for the Hill (black data) and the autocatalytic (blue data) models using log-log scale and $x(0) = 0.75\cdot \Omega$. Here, $\phi^{(s)} = \epsilon - \bar{\epsilon}^{(s)}$ is the distance between the bifurcation parameter, $\epsilon$, and the mean stochastic bifurcation value, $\bar{\epsilon}^{(s)}$, see~\cite{SardanyesNJP2020}. Close to the stochastic bifurcation transients' scaling remains flat (red rectangle). (Bottom) Scaling obtained for the flight times computed from the Hamiltonian system, averaging over $10^2$ different initial conditions (variances are not shown for the sake of clarity) within the region of lowest action. In all of the analyses we used $\Omega = 10^3$. We show different representative dynamical regimes: (a) long transients for $\phi =1.127 \times 10^{-5}$, (b) strong change of the curvature of the flat region at $\phi \approx 0.0011$, and shorter transients with (c) $\phi =10^{-2}$ and (d) $\phi  \approx 0.114$.}
\label{fig_1}
\end{center}
\end{figure*}
%FIGURE %%%%%%%%%%%%%%%%%%%%%%%%%%%%%%%%%%%%%%%%%%%%%%%%%%%%%%%%%%%%

\subsection{Scaling of the extinction times close to the critical point}
\label{scaling:Hill}
We have first computed the scaling of the mean extinction times, $\bar{T}_E$, with GS for the two representative models under consideration (see Table \ref{tab:table_1} and~\cite{SardanyesNJP2020}). Specifically, we have studied how $\bar{T}_E$ changes as $\phi^{(s)} = \epsilon - \bar{\epsilon}^{(s)}$ varies. That is, as the control parameter is driven beyond the mean stochastic bifurcation value, which has been computed using the same approach as in~\cite{SardanyesNJP2020}. The first observation is that, close to the bifurcation, the behaviour of both models is characterised by the same scaling function (see upper panel in Fig.~\ref{fig_1}). As expected, far from the critical point the scaling behaviour is lost, and $\bar{T}_E$ shows model-specific features.

\bigskip
%FIGURE %%%%%%%%%%%%%%%%%%%%%%%%%%%%%%%%%%%%%%%%%%%%%%%%%%%%%%%%%%%%
\begin{figure*}
\begin{center}
\centerline{\includegraphics[width=0.925\textwidth]{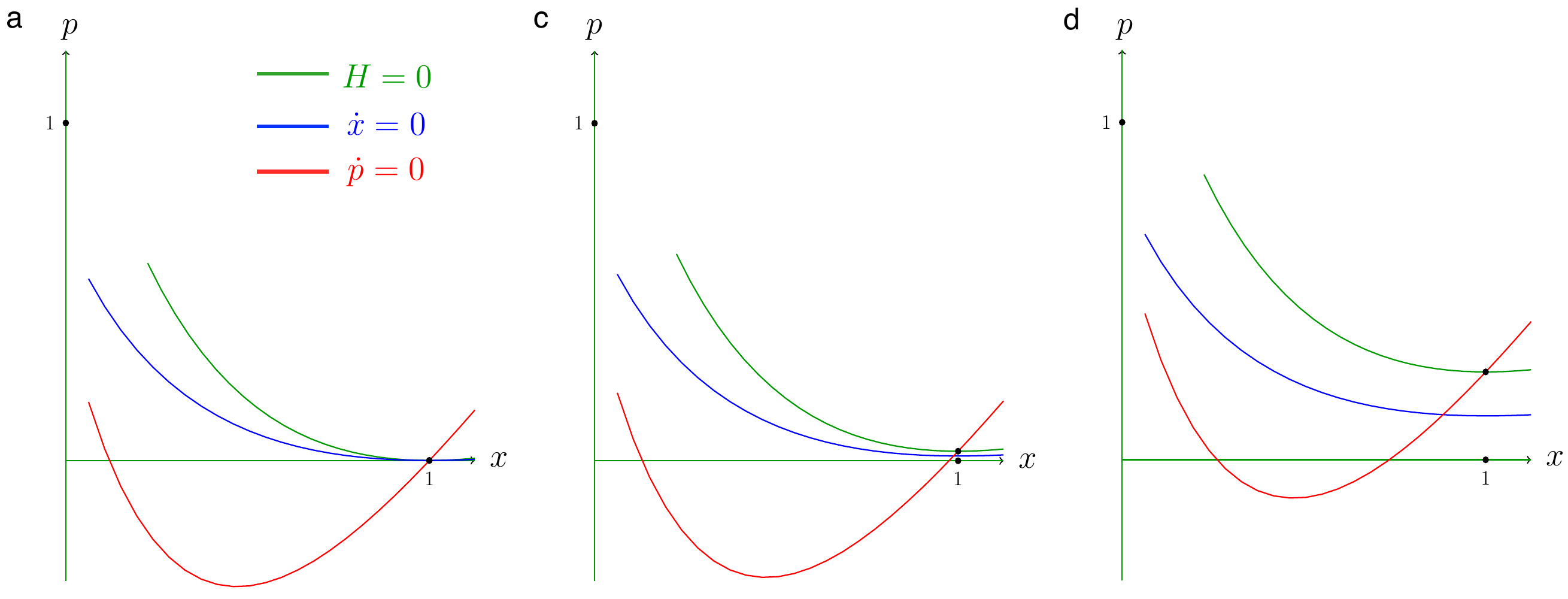}}
\caption{Relevant curves determining the dynamics of the Hamiltonian Hill model given by Eqs.~\eqref{eq:Hill:A:1}. The nullclines $\dot{x}=0$ (blue) and  $\dot{p}=0$ (red) are shown together with the zero value of the Hamiltonian curve (green). The three plots follow the same labels as in Fig. 1: (a) $\phi =1.127 \times 10^{-5}$, (c) $\phi=10^{-2}$, and (d) $\phi  \approx 0.114$.}
\label{fig_2}
\end{center}
\end{figure*}
%FIGURE %%%%%%%%%%%%%%%%%%%%%%%%%%%%%%%%%%%%%%%%%%%%%%%%%%%%%%%%%%%%

The scaling obtained with the GS is well captured by the semiclassical theory (see lower panel in Fig.~\ref{fig_1}). We have solved numerically\footnote{Numerical integration has been carried out with a Runge-Kutta-Fehlberg method of order 7-8 with adaptive time step size and local tolerance $10^{-15}$.} the Hamilton equations~\eqref{eq:hamiltoneqsp} for the two studied systems (shown in Sections \ref{A11} and \ref{se:autocatalytic:model}) in order to compute the extinction times associated with their orbits for different values of $\phi$. Our results show that the same shape, i.e. the same scaling function, is obtained. The scaling for the flight times, which do not depend on $\Omega$, at changing $\phi$ displayed in Fig.~\ref{fig_1} match with the GS ran  close to bifurcation threshold (red rectangle) using $\Omega = 10^3$. Notice that the flat shape close to this threshold starts bending at $\phi \approx 0.00043$ for this system size. We note that using other system sizes also drive to good matching between the GS and the flight times obtained numerically from the Hamiltonians. %\textcolor{blue}{Se puede escribir de dos maneras: (i) For instance, the change from the flat region towards the bending of the curve occurs at $\phi \approx 0.0008$ for $\Omega = 8\times 10^3$; $\phi \approx 0.0006$ for $\Omega = 4\times 10^3$; $\phi \approx 0.0005$ for $\Omega = 2\times 10^3$; and $\phi \approx 0.0003$ for $\Omega = 500$ (results not shown, see also \cite{SardanyesNJP2020}).}
For instance, system sizes ranging from  $\Omega = 500$ to $\Omega = 8\times 10^3$ show, respectively, values of $\phi$ causing the bending of the flat region within the approximate range $\phi \in [0.00035, 0.00072]$ (results not shown, see also \cite{SardanyesNJP2020}).

%To gain some further understanding on these observations, we have solved numerically the Hamilton equations of the function $H(x,p)$ (see Eqs.~\eqref{eq:hamiltoneqsp}). 
%% TL (to check)
The two investigated models, the Hill and the autocatalytic, share a phase space structure which clearly determines their dynamics (see Section \ref{A1} for technical details)). For example, for the Hill model, such configuration is determined by the main curves: $H(x,p)=0$ and the two nullclines $\dot{x}=0$ and $\dot{p}=0$. The first one, is constituted by three components, the coordinate axes and the curve $y=p_H(x)=\log \left( \eps \frac{1+x^2}{x} \right)$. This curve is tangent to $p=0$ exactly at the s-n bifurcation value $\eps=\eps_c$, when it meets both nullclines at $x=x_c$ (see Figure~\ref{fig_2}a).  
As we increase the value of $\phi$, the curve $y=p_H(x)$ moves upwards and gives rise to a 'tunnel' of width $\mathcal{O}(\phi)$ which embraces the $x$-nullcline - given by $p=p_1(x)=\frac{1}{2} p_H(x)$ (see Figure~\ref{fig_2}b). The $p$-nullclines $p=0$ and $p=p_2(x)=\log \left( \frac{\eps (1+x^2)^2}{2x} \right)$ and the regions of growth/decline which all them ($\dot{x}=0$, $\dot{p}=0$ and $H=0$) force the dynamics for $p\gtrsim 0$ to evolve $\mathcal{O}(\phi)$-close to the line $p=0$ (deterministic case). This is not the case for orbits with $p\lesssim 0$, for which the dynamics inside the lobe determined by $p=0$ and $p=p_2(x)$ can exhibit a $\mathcal{O}(1)$-evolution in $p$, small but larger than the $\mathcal{O}(\phi)$ observed above (see Figure~\ref{fig_2}c). This configuration disappears once the curve $p=p_2(x)$, which moves also upwards as a function of $\phi$, presents a tangency with $p=0$ for a value of $\phi$ of order $5 \cdot 10^{-2}$ approximately (see Figure~\ref{fig_1}, bottom).     

It is worth mentioning that, close to the bifurcation, the behaviour of the extinction time depends crucially on the initial condition of the momentum variable, $p_0$. If $p_0\gtrsim 0$, then the extinction time follows the deterministic scaling law (i.e. $T_E\sim \phi^{-1/2}$, results not shown). By contrast, orbits with $p_0\lesssim 0$ exhibit a scaling function whose shape is the same as the GS.  This difference in behaviour appears to be predicated upon the feature of the phase space $(x, p)$ shown in Fig. \ref{fig_2}. When $\phi\rightarrow 0$ and $p_0\lesssim 0$, the orbits of the system cross the region of the phase space contained within the red line representing the nullcline $\dot{p}=0$. This sojourn alters the properties of the slowing down, thus changing the scaling behaviour of the transit time. Section~\ref{scalings_ps} provides analytical results explaining the changes in the scaling functions for positive and negative $p_0$.
%FIGURE %%%%%%%%%%%%%%%%%%%%%%%%%%%%%%%%%%%%%%%%%%%%%%%%%%%%%%%%%%%%
\begin{figure*}
\begin{center}
\centerline{\includegraphics[width=0.8\textwidth]{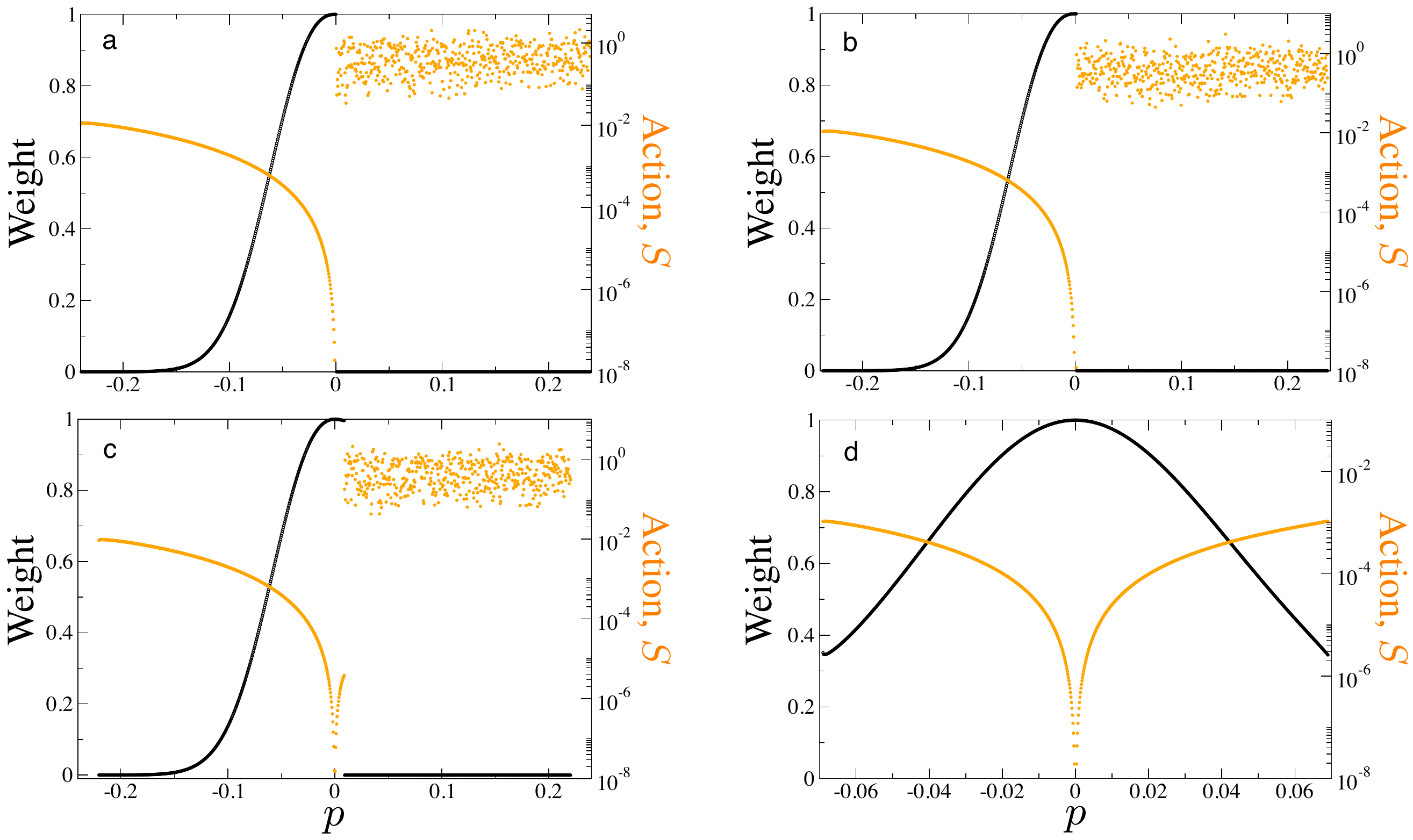}}
\caption{Weight of the regions in $p$ along which stochastic realisations will pass (black lines) for the Hill model, and  the associated action to the values of variable $p\equiv p(t=0)$ (orange dots). The four panels show the system's behaviour for the values of $\phi$ indicated with the same letters in Fig. 1. Here we have used $\Omega = 10^3$.}
\label{fig_3}
\end{center}
\end{figure*}
%FIGURE %%%%%%%%%%%%%%%%%%%%%%%%%%%%%%%%%%%%%%%%%%%%%%%%%%%%%%%%%%%%

\subsection{Statistical significance of the orbits of the semiclassical Hamiltonian}

In order to shed some light onto the picture described in Section~\ref{scaling:Hill},  we have analysed the statistical significance of the different sets of orbits of $H(x,p)$ for the Hill model. According to the path integral representation of the stochastic process, the weight of each trajectory is determined by their associated action, $S$ [see Eq.~\eqref{eq:hamiltonian}]: $e^{-\Omega S}$. We show that close to the bifurcation the only statistically significant trajectories are the ones with initial $p_0\lesssim 0$, as shown in Fig.~\ref{fig_3}(a,b). These paths are the ones that contribute to the stochastic scaling function. By contrast, as we move away from the critical point, both orbits with positive and negative initial values of $p$ become statistically significant, as shown in Fig.~\ref{fig_3}(c) with low action (weight $\sim 1$) for small, positive $p$. Eventually [see Fig.~\ref{fig_3}(d)], trajectories with both positive and negative $p_0$ have equal weights and therefore the behaviour of the system becomes a mixture of both characteristic behaviours.

\section{Discussion}

In this paper we have analysed the dynamical behaviour of stochastic systems close to a first-order phase transition, mainly regarding the slowing down phenomena observed in the relaxation times towards their absorbing state. We recently explored the effect of intrinsic noise in simple models with cooperation undergoing s-n bifurcations including the autocatalytic replicator and two-species hypercycles~\cite{SardanyesNJP2020}. Somewhat counter-intuitively, we found that the deterministic ghosts were robust to noise. However, the scaling law of the mean extinction times was shown to follow a rather more complicated scaling function than its deterministic counterpart [see Eq.~\eqref{eq:TE}]. 

In an attempt to better understand the mechanisms underlying this behaviour, in particular concerning the origin of the scaling function, ${\mathcal{G}}(\cdot)$, we have reformulated the Master Equation (ME) as a path integral. Hence the solution of the ME is obtained as a weighted superposition of paths in the phase space of the Hamiltonian (\ref{eq:hamiltonian}), which are given by the solutions of Eqs.~\eqref{eq:hamiltoneqsp}. Each such path is assigned a weight that depends exponentially on the action associated with it [see Eq. \eqref{eq:pathint}] times the system size $\Omega$.

The analysis of the paths and their weights yields a number of interesting results which shed some light onto the origin of both the robustness of the ghost behaviour in the presence of intrinsic noise and the scaling function ${\mathcal{G}}$. First, when the system is close to the bifurcation ($\phi=\epsilon-\epsilon_c={\mathcal{O}}(10^{-3})$ or smaller), there is a set of solutions (corresponding to $p(t=0)<0$) of the Hamiltonian system, Eqs.~\eqref{eq:hamiltoneqsp}, which undergo hugely long transients when trying to cross the regions of the phase space shown in Fig. \ref{fig_2}. Furthermore, as shown in Fig. \ref{fig_1} (lower panel), when we compute how the flight time associated with the paths varies as $\phi$ changes, we realise that, rather than following the deterministic scaling law $\phi^{-1/2}$, it is qualitatively the same as the one computed numerically by means of Gillespie simulations for values of $\phi\leq 10^{-3}$ (see Fig. \ref{fig_1}, upper panel). By contrast, solutions of Eqs.~\eqref{eq:hamiltoneqsp} with $p(t=0)\geq 0$ exhibit the same scaling behaviour than the deterministic limit, i.e. flight time $\sim\phi^{-1/2}$ (results not shown, see also Section~\ref{scalings_ps}). 
Further insight is obtained from the quantification of the action and weight associated with each path as shown in Fig. \ref{fig_3}. These results show that, for $\phi\leq 10^{-3}$, the only statistically significant paths are those corresponding to $p(t=0)<0$, whereas all the paths with $p(t=0)\geq 0$ have negligibly small weights. Taken together, these results suggest that, close to the bifurcation, the only statistically significant paths, and consequently the only ones that will contribute to the observable behaviour, are those with $p(t=0)<0$, for which the flight time is of the order a (small) constant.
%$\simeq{\mathcal{O}}(\phi)$. 

Our core outcome is therefore that the qualitative behaviours displayed by the paths of the Hamiltonian Eq.~\eqref{eq:hamiltonian}, are thus mimicked by or translated to the observable behaviour of the stochastic system, provided that they are statistically significant. In the specific case of systems whose mean-field model exhibits a s-n bifurcation, the delayed relaxation and shape of the function ${\mathcal{G}}$ are the statistically significant features that are passed on to the stochastic system. This explains the robustness of the ghost behaviour to intrinsic noise and the scaling behaviour of the mean extinction times, $\bar{T}_E$ found in~\cite{SardanyesNJP2020}. It is noteworthy that, although we have illustrated our theory with two specific examples, the same properties should extend in a universal way to all stochastic systems whose associated Hamiltonian system exhibits the same behaviour.    

Finally, we conclude that the analysis we have put forward here can have great relevance to analyse qualitative behaviours of complex stochastic systems, specifically, close to phase transitions/bifurcations, in a way that does not require a substantial computational effort. By analysing the orbits of the corresponding Hamiltonian system and the different qualitative behaviours it can sustain, we will be able to evaluate which behaviours are associated to statistically relevant orbits, and therefore should manifest themselves at the level of the observable behaviour of the stochastic system.

%In This scaling relation was identified in a simple autocatalytic model and in a two-species catalytic network~\cite{SardanyesNJP2020}. These findings were obtained by means of extensive numerical simulations and thus a theory explaining this scaling relations remained undeveloped. We here address this issue using the WKB asymptotic approximation of the Master Equation, which allows to provide a dynamical investigation of this phenomenon in the phase space of the Hamiltonian obtained from the Hamilton-Jacobi equation. By applying this approach to the autocatalytic replicator model studied in~\cite{SardanyesNJP2020} as well as to a simple model of cooperative gene regulation, we have found scalings function matching those obtained from the stochastic numerical simulations for such systems.  

\section*{Acknowledgements} 
 TA, JS and TL would like to thank the Laboratorio Subterr\'aneo de Canfranc (LSC) for kind hospitality during the development of this research. This work is supported by the Spanish State Research Agency (AEI), through the Severo Ochoa and María de Maeztu Program for Centers and Units of Excellence in R\&D (CEX2020-001084-M). We thank CERCA Programme/Generalitat de Catalunya for institutional support. JS and TA acknowledge the AEI for funding under grant RTI2018-098322-B-I00. JS has been also funded by the Ram\'on y Cajal contract RYC-2017-22243 and TL by the Spanish project PGC2018-098676-B-100.

%\section*{References}

\appendix

\section{Analytic results}

\subsection{Qualitative features of the Hamiltonian system} \label{A1}

Here we discuss the qualitative features of the semiclassical Hamiltonian systems associated with systems exhibiting a saddle-node bifurcation. Specifically, we analyse in detail the properties of the nullclines shown in Fig. \ref{fig_2}.

\subsubsection{Hill model.} 
\label{A11}
According to Table \ref{tab:table_1} and to Eqs. (\ref{eq:scalingrates}), (\ref{eq:hamiltonian}), (\ref{eq:hamiltoneqsp}), and performing the change of variables  $(x,p)\mapsto (y=\frac{x}{A},p)$, time $t= A \sigma$, as well as  redefining the parameter $\eps\rightarrow A \eps$, we have
\begin{equation}
\begin{array}{rcl}
\dot{x} &=& {\ds - \rme^p \frac{x^2}{1 +x^2} - \rme^{-p} \eps x} , \label{eq:Hill:A:1} \\
\dot{p} &=& {\ds - \left( \rme^p - 1 \right) \, \frac{2x}{\left( 1+x^2\right)^2} - \left( \rme^{-p} -1 \right) \eps}.
\end{array}
\end{equation}
The deterministic limit (i.e., for $p=0$) has three equilibrium points: $x=0$ (attracting) and 
$\ds x_{\pm} = \frac{1}{2\eps}\left ( 1 \pm  \sqrt{1 -4\eps^2}\right).$
The saddle-node bifurcation occurs at $\eps= \criteps=\frac{1}{2}$. At the critical point, the non-trivial equilibrium is
$x_{\mathrm{c}}= \frac{1}{2\criteps}= 1.$ The relevant curves (see Fig. \ref{fig_2}) are:
\begin{itemize}
\item The solution of $H(x,y)=0$ is formed by three curves, namely, $x=0$, $p=0$, and $p=p_H(x)$, where
$$
p_H(x):=\log \left( \eps \, \frac{1+x^2}{x} \right).  
$$
Note that $y=p_H(x)$ intersects the axis $p=0$ only for $\epsilon=\criteps$. In fact, $p=p_H(x)$ goes to $+\infty$ as $x\rightarrow 0^+$  and as $x\rightarrow +\infty$, and it has a unique global minimum,  $x_{\textrm{min,H}}=1$ with $p_{\textrm{min,H}}= p_{H}(x_{\textrm{min,H}}) = \log 2\epsilon.$

\item The $x$-nullcline $\dot{x}=0$ is formed by the curves $x=0$ and 
$p=p_1(x):= \frac{1}{2} p_H(x)$, 
%\label{eq:x-nullcline}
which exhibit a minimum value 
$p_1(x_{\textrm{min,H}})= \frac{1}{2} \log(2\eps)$.

\item The $p$-nullcline $\dot{p}=0$ has two components: $p=0$ and $p=p_2(x)$, with
$$
p_2(x):=\log \left( \frac{\eps (1+x^2)^2}{2x} \right)  = p_H(x) + \log \left( \frac{1+x^2}{2} \right).
$$
The curve $p=p_2(x)$ has a global minimum at
$x_{\textrm{min,p}}=\frac{1}{\sqrt{3}}$ with 
$p_{\textrm{min,p}} = p_2(x_{\textrm{min,p}})=\log \left( \frac{8\sqrt{3}}{9}\eps \right)$.
The intersection points $x_F(\eps) \leq x_0(\eps)$ of $p=p_2(x)$ with $p=0$ are given by the two (positive) solutions of $g(x)=\eps$, where
$g(x) = \frac{2x}{\left( 1+x^2\right)^2}.$ One can show that the maximal interval of $\eps$ where $g(x)=\eps$ has feasible solutions is
\begin{equation}
\eps \in [\criteps, \endeps] := \left[ \frac{1}{2} \, , \, \frac{3\sqrt{3}}{8} \right], 
\label{eq:epsilon:interval}
\end{equation}
so that as long as $\eps$ stays within this interval the system can exhibit the scaling behaviour shown in Fig. \ref{fig_1}.
\end{itemize}

\subsubsection{Autocatalytic model.}
\label{se:autocatalytic:model}

Similarly, for the autocatalytic model, the relevant curves are:
\begin{itemize}
\item The solution of $H=0$ is constituted by three curves, namely, $x=0$, $p=0$ and $p=p_H(x)$, where
$$
p_H(x):=\log \left( \frac{\eps + x^2}{x} \right)
$$
has a global minimum at 
$x_{\textrm{min,H}}=\sqrt{\eps},$
$p_{\textrm{min,H}}=p_H(x_{\textrm{min,H}}) = \log(2\sqrt{\eps}).$
The value of $\eps$ for which $p=p_H(x)$ reaches its minimum on $p=0$ is $\eps=\criteps=\frac{1}{4}$, and the point is $x_{\textrm{min,H}} = x_{\textrm{c}}=\frac{1}{2}$. 

\item The $x$-nullcline, $\dot{x}=0$, has two branches: the curve $x=0$ and
$p=p_1(x)$, where
$$
p_1(x):=\frac{1}{2} \log \left( \frac{x^2+\eps}{x} \right) = \frac{1}{2} p_H(x). 
$$

\item The $p$-nullcline, $\dot{p}=0$, also has two branches: $p=0$ and
$p=p_2(x)$, where
$$
p_2(x):=\log \left( \frac{3x^2+\eps}{2x} \right), 
$$
which has a global minimum at
$x_{\textrm{min,p}}=\sqrt{\frac{\eps}{3}},$ and 
$p_{\textrm{min,p}}=p_2(x_{\textrm{min,p}})= \log\sqrt{3\eps}.$
Moreover, $p=p_2(x)$ is tangent to $p=0$ for $\eps=\endeps=\frac{1}{3}$. For $\eps>\endeps$ one has $p_2(x)>0$. This leads to
\begin{equation}
\eps \in [\criteps, \endeps] := \left[ \frac{1}{4} \, , \frac{1}{3} \right] 
\label{eq:epsilon:interval:autocat}
\end{equation}
as the interval in the bifurcation parameter $\eps$ where the slowing down is studied.
\end{itemize}

\subsection{On the scaling law for $p>0$ and $p<0$} \label{scalings_ps}
Let us restrict ourselves to the autocatalytic model. A similar approach can be performed for the Hill equation. We essentially follow the ideas presented in~\cite{Fontich2008}.

From the first equation,
$$
\dot{x}=x^2 \rme^p - \left( x^3 + \eps x \right) \rme^{-p} =: \tilde{f}(x,p,\eps)
$$
it follows that the time $T(\eps)$ invested by a solution with initial condition on $x=x_{\textrm{ini}}$ to reach $x=x_{\textrm{end}}$ is given by
$$
T(\eps) = \int_{x_{\textrm{ini}}}^{x_{\textrm{end}}} \frac{1}{\tilde{f}(x,p,\eps)} \, dx.
$$
Our aim is to compute its first order approximation in $\phi=\eps -\criteps$. To do it, first, notice that
$$
x^2 \rme^p - \left( x^3 + \eps x \right) \rme^{-p} = -x \rme^{-p} \left( x^2 - \rme^{2p} x + \eps \right) = \left( x-\frac{1}{2} \right)^2 + \left( 1-\rme^{2p}\right) x + \phi, 
$$
where $\phi = \eps - \frac{1}{4}$, being $\eps_{\textrm{c}}=\frac{1}{4}$ the value of $\eps$ at which the systems undergoes a saddle-node bifurcation.
Remind that the corresponding $x$-value for $\eps=\eps_{\textrm{c}}$ is $x_{\textrm{c}}=\frac{1}{2}$.
If we apply the change of variables $y=x - \frac{1}{2}$ then $y(t)$ satisfies the following ODE:
\begin{equation}
\dot{y}=-\rme^{-p} \left( y + \frac{1}{2} \right) \bigg( y^2 + \left( 1-\rme^{2p} \right) \left( y + \frac{1}{2} \right) + \phi \bigg) =: f(y,p,\phi). 
\label{ode:y}
\end{equation}
Observe that
\begin{equation}
f(0,0,0)=0 \qquad \textrm{and} \qquad \frac{\partial f}{\partial \phi}(y,p,\phi)=-\rme^{-p} \left( y + \frac{1}{2} \right) \Rightarrow \frac{\partial f}{\partial \phi}(0,0,0)=-\frac{1}{2} \ne 0.   
\label{IFT:conditions}
\end{equation}
\begin{lemma}
There exists an open set $\mathcal{U} \in \R^2$, containing $(0,0)$, and a $\mathcal{C}^{\infty}$-function 
$$
\begin{array}{rcll}
g : &\mathcal{U} & \longrightarrow & \R \\
    & (y,p) & \longmapsto & \phi=g(y,p)  
\end{array}
$$
such that $g:\mathcal{U} \rightarrow g(\mathcal{U})$ is bijective and satisfies that $f(y,p,\phi)=0 \Leftrightarrow \phi=g(y,p)$ $\forall (y,p)\in \mathcal{U}$. 

Let us denote, by commodity, $\phi=\phi(y,p)=g(y,p)$.
Then, $\phi(y,p)$ admits the following Taylor expansion around $(0,0)$:
\begin{equation}
\phi(y,p)=p - y^2 + 2py + p^2 + \mathcal{O}_3(y,p),
\label{phi:taylor}
\end{equation}
where $\mathcal{O}_3(y,p)$ stands for terms of order, at least $3$, in $y$ and $p$. 
\end{lemma}

\bigskip

\noindent\textbf{Proof.}
Conditions~\eqref{IFT:conditions} allow us to apply the Implicit Function Theorem and prove the first assertion of the lemma. Regarding the second one, we seek for
\begin{eqnarray}
\phi(y,p) &=&\phi(0,0) + \frac{\partial \phi}{\partial y}(0,0)\, y + 
\frac{\partial \phi}{\partial p}(0,0)\, p \nonumber \\
&+& \frac{1}{2} \left( 
\frac{\partial^2 \phi}{\partial y^2}(0,0)\, y^2 +
2 \frac{\partial^2 \phi}{\partial y \partial p}(0,0)\, yp +
\frac{\partial^2 \phi}{\partial p^2}(0,0)\, p^2 \right) + \mathcal{O}_3(y,p). 
\label{taylor:1}
\end{eqnarray}
Standard recurrent differentiation of equation $f(y,p,\phi)=0$ - with respect to $y$ and $p$ - and substitution onto $y=p=0$ gives $\phi_y(0,0)=0$, $\phi_p(0,0)=0$, $\phi_{yp}(0,0)=2$, $\phi_{yy}(0,0)=-2$, and $\phi_{pp}(0,0)=2$. This leads to expression~\eqref{phi:taylor}. 

\qed

\bigskip

From
$\phi=p - y^2 + 2py + p^2 + \mathcal{O}_3(y,p)$,
it follows that
$0= y^2 - 2py-p^2 -p + \phi + \mathcal{O}(\phi^3) \Rightarrow 
(y-p)^2 = 2p^2 + p - \phi + \mathcal{O}(\phi^3)$
and so, 
$$
f(y,p,\phi)=0 \Leftrightarrow y= y_{\pm}= p \pm \sqrt{2p^2+p-\phi} + \mathcal{O}(\phi). 
$$
Notice that for $p=0$ one has~\footnote{Sometimes, to emphasize its dependence on $\phi$ we will write $y_{\pm}(\phi)$.} $y_{\pm}=\pm \rmi \sqrt{\phi} + \mathcal{O}(\phi)$. Hence, for $p\sim 0$ we have that $y_{\pm}$ are close to
$y_{\pm}^0 = \pm \rmi \sqrt{\phi}$. We consider the case when $2p^2+p-\phi$ is negative and so $y_{\pm} \sim y_{\pm}^0$ are complex. Hence
$
2p^2+p-\phi = 2 (p-p_{-}) (p-p_{+}), 
$
where
$$
p_{\pm}=\frac{-1\pm \sqrt{1+8\phi}}{4} \in \R \qquad \textrm{since $\phi>0$}. 
$$
Precisely $p_{-} < -\frac{1}{4}$ and
\begin{eqnarray*}
p_{+}&=&\frac{-1 + \sqrt{1+8\phi}}{4} =
\frac{-1 + \sqrt{1+8\phi}}{4} \cdot \frac{1 + \sqrt{1+8\phi}}{1 + \sqrt{1+8\phi}} = \frac{2\phi}{1 + \sqrt{1+8\phi}} \\
&=& \frac{2\phi}{1+\left( 1 - 4\phi + \mathcal{O}(\phi^2) \right)} =
\frac{\phi}{1 - 2\phi + \mathcal{O}(\phi^2)} \\[1.2ex]
&=& \phi \left( 1+2\phi + \mathcal{O}(\phi^2) \right) = \phi + \mathcal{O}(\phi^2). 
\end{eqnarray*}
Roughtly speaking, $p_{-}=\mathcal{O}(1)$ and $p_{+}=\mathcal{O}(\phi)>0$. Therefore,
$$
y_{\pm} = p \pm 2\rmi \sqrt{p-p_{-}} \, \sqrt{p-p_{+}} + \mathcal{O}(\phi),
$$
for $p_{-} < p < p_{+}$. Moreover, $\Im y_{+}>0$ and $\Im y_{-}<0.$
We discuss the cases $p>0$ and $p<0$ separately.

\bigskip

\textbf{Case $\mathbf{p>0}$.} As mentioned at the beginning of this section, our aim is to estimate the flight time of solutions crossing close to the section $x=x_{\textrm{c}}=\frac{1}{2}$. 

Let us take $\phi_0 \in (0,\frac{1}{12}]$, where $\frac{1}{12}=\eps_{ \textrm{end}}-\eps_{ \textrm{c}}=\frac{1}{3} - \frac{1}{4}$, the $\phi$-interval for which the curve $\dot{p}=0$ intersects the line $p=0$.
Now fix $\phi\in (0,\phi_0]$. From the change of variables $y=x-\frac{1}{2}$ and having in mind that our system is autonomous, the flight time $T(\phi)$ of trajectories starting at $y=\delta$ and reaching $y=-\delta$ can be computed as
\begin{equation}
T_+(\phi) = \int_{\delta}^{-\delta} \frac{1}{f(y,p,\phi)} \, dy,
\end{equation}
for some $\delta>0$ small. Let consider $0<\nu<p_+(\phi_0)$, independent of $\phi$, such that  
for $(y,p) \in W=[-\delta,\delta] \times (0,\rmi \nu]$ the function $f(y,p,\phi)$ does not vanish, except at the point $y=y_{+}(\phi) \simeq \rmi\sqrt{\phi}$. In other words, $\frac{1}{f}$ is meromorphic in $W$ and has a unique pole $y_{+}\in \mathring{W}$. Furthermore, this pole is simple. 
Let denote by $\gamma=\partial^+W$ the (positively oriented) border of $W$. Thus, $\gamma=\gamma_1 \cup \gamma_2 \cup \gamma_3 \cup \gamma_4$, being $\gamma_1$ the segment placed on $p=0$, $\gamma_2$ the one on $y=\delta$, $\gamma_3$ the one on $p=\rmi\nu$, and $\gamma_4$ the one on $y=-\delta$. Therefore, 
\begin{equation}
T_{+}(\phi) = -\int_{\gamma_1} \frac{1}{f(y,p,\phi)} \, dy = -\int_{\gamma} \frac{1}{f(y,p,\phi)} \, dy + \sum_{j=2}^4 \int_{\gamma_j} \frac{1}{f(y,p,\phi)} \, dy. 
\label{p:positive:time:1}
\end{equation}
By the Residue Theorem we have that
\begin{equation}
\int_{\gamma} \frac{1}{f(y,p,\phi)} \, dy = 2\pi\rmi \, \mathrm{Res}\left( \frac{1}{f}, y=y_{+}\right) = 2\pi\rmi \lim_{y\rightarrow y_{+}} \frac{y-y_{+}}{f(y,p,\phi)} = \frac{2\pi\rmi}{\frac{\partial f}{\partial y} (y_{+},p,\phi)}.
\label{p:positive:residue}
\end{equation}
Notice that
$$
y_{+} = p + 2\rmi \sqrt{p-p_{-}}\, \sqrt{p_{+}-p} = 2\rmi c_1 \sqrt{\phi} + \mathcal{O}(\phi), 
$$
since $p_+=\mathcal{O}(\phi)$, $p=\mathcal{O}(\phi)$, $\sqrt{p_{+}-p}=\mathcal{O}(\sqrt{\phi})$, and $\sqrt{p-p_{-}}=c_1=\mathcal{O}(1)$. Substituting the expression $y_{+}=2\rmi c_1 \sqrt{\phi} + \mathcal{O}(\phi)$ into
\begin{equation}
\frac{\partial f}{\partial y}(y,p,\phi) = -\rme^{-p} \left( 3y^2 + 2 \left( 1-\rme^{2p} \right) \left(y+\frac{1}{2} \right) + y + \phi \right) 
\label{dfdy}
\end{equation}
leads to 
$$
\frac{\partial f}{\partial y}(y_{+},p,\phi) = -2\rmi c_1 \sqrt{\phi} + \mathcal{O}(\phi) 
$$
and so equation~\eqref{p:positive:residue} becomes
\begin{equation}
\int_{\gamma} \frac{1}{f(y,p,\phi)} \, dy =
\frac{1}{\frac{\partial f}{\partial y} (y_{+},p,\phi)} = \frac{2\pi\rmi}{-2\rmi c_1 \sqrt{\phi} + \mathcal{O}(\phi)} = -\frac{\pi}{c_1} \, \frac{1}{\sqrt{\phi}} + \mathcal{O}(\phi). 
\label{p:positive:int:gamma1}
\end{equation}
Let us now deal with the computation on the curves $\gamma_2, \gamma_3$, and $\gamma_4$.  We start with $\gamma_2$, which can be parameterisated by $\gamma_2(\sigma)=\delta + \rmi \nu \sigma$ for $\sigma\in [0,1]$. Since
$f(\gamma_2(\sigma),\phi)=f(\gamma_2(\sigma),0)+\mathcal{O}(\phi)$ and by
the compactness of the interval $[0,1]$ it follows that
$$
\frac{1}{f(\gamma_2(\sigma),\phi)} = \frac{1}{f(\gamma_2(\sigma),0)} + \mathcal{O}(\phi).  
$$
Therefore,
\begin{equation}
\int_{\gamma_2} \frac{1}{f(y,p,\phi)} \, dy =
\int_0^1 \frac{1}{f(\gamma_2(\sigma),\phi)} \, d\sigma =
\int_0^1 \left( \frac{1}{f(\gamma_2(\sigma),0)} + \mathcal{O}(\phi) \right) \, d\sigma = c_2 + \mathcal{O}(\phi),
\label{p:positive:gamma2}
\end{equation}
$c_2$ being a constant. The same argument (and similar result) applies for the integrals over $\gamma_3$ and $\gamma_4$. In those cases, the parameterisations being 
$\gamma_3(\sigma)=(2\sigma-1) \delta + \rmi \nu$ and
$\gamma_4(\sigma)=-\delta + \rmi \nu (1-\sigma)$, with $\sigma \in [0,1]$, respectively. Having these estimates in mind and expressions~\eqref{p:positive:time:1} and~\eqref{p:positive:int:gamma1} we obtain that
\begin{equation}
T_{+}(\phi) = \frac{\pi}{c_1} \, \frac{1}{\sqrt{\phi}} + C +\mathcal{O}(\phi),  
\end{equation}
where $C=c_2+c_3+c_4$, a constant independent of $\phi$. In logarithmic scale this relation reads as
$$
\log T_+(\phi) \sim -\frac{1}{2} \log \phi, 
$$
the expected scaling law for a saddle-node bifurcation in a deterministic system (see~\cite{Fontich2008}).

\bigskip

\noindent\textbf{Case $\mathbf{p<0}$.} As for $p>0$, take $\phi_0 \in (0,\frac{1}{12}]$ and fix $\phi \in (0,\phi_0]$. To compute
$$
T_{-}(\phi) = \int_{\delta}^{-\delta} \frac{1}{f(y,p,\phi)} \, dy 
$$
we consider $0<\eta<\frac{1}{4}$, independent of $\phi$ and $\eta = \mathcal{O}(1)$, such that $f(y,p,\phi)\ne 0$ for any $(y,p)\in V=[-\delta,\delta] \times [-\rmi\eta,0)$ except the point $y_{-}(\phi)$, where it has a simple zero. As above, this means that $\frac{1}{f}$ is meromorphic in $V$ and has a unique pole (simple) at $y=y_{-}(\phi)\in \mathring{V}$. We follow the same argument as in the previous case: define $\Gamma=\Gamma_1\cup \Gamma_2 \cup \Gamma_3 \Gamma_4$, the positively oriented border of $V$. $\Gamma_1$ is the segment on $p=0$, $\Gamma_2$ the one on $y=-\delta$, $\Gamma_3$ the one on $p=-\rmi\eta$, and $\Gamma_4$ the one on $y=\delta$. Therefore, 
\begin{equation}
T_{-}(\phi) = \int_{\Gamma_1} \frac{1}{f(y,p,\phi)} \, dy = \int_{\Gamma} \frac{1}{f(y,p,\phi)} \, dy - \sum_{j=2}^4 \int_{\Gamma_j} \frac{1}{f(y,p,\phi)} \, dy. 
\label{p:negative:time:1}
\end{equation}
Like in the previous case, we compute the first line integral by the Residue Theorem:
\begin{equation}
\int_{\Gamma} \frac{1}{f(y,p,\phi)} \, dy = 2\pi\rmi \, \mathrm{Res}\left( \frac{1}{f}, y=y_{-}\right) = 2\pi\rmi \lim_{y\rightarrow y_{-}} \frac{y-y_{-}}{f(y,p,\phi)} = \frac{2\pi\rmi}{\frac{\partial f}{\partial y} (y_{-},p,\phi)}.
\label{p:negative:residue}
\end{equation}
To get an estimate for this value we use the definition for 
$\frac{\partial f}{\partial y} (y_{-},p,\phi)$ given in Equation~\ref{dfdy}. Contrary to what occurred for $p>0$, restricted to move in a strip of order $\mathcal{O}(\phi)$, now $p<0$ evolve at a distance of order $\mathcal{O}(1)$ from $p_{-}<-\frac{1}{4}$ and from $p_+>0$. This makes that $p=\mathcal{O}(1)$, as well as $y_{-}$, $\rme^{p}$, and $2(1-\rme^{2p})$. Hence
$$
\frac{2\pi\rmi}{\frac{\partial f}{\partial y} (y_{-},p,\phi)} =
\frac{- 2\pi\rmi \rme^p}{3y_{-}^2 + 2 \left( 1-\rme^{2p} \right) \left(y_{-}+\frac{1}{2} \right) + y_{-} + \phi} = \mathcal{O}(1) + \mathcal{O}(\phi).
$$
Consequently, 
\begin{equation}
\int_{\Gamma} \frac{1}{f(y,p,\phi)} \, dy =  c_5 +  \mathcal{O}(\phi), 
\label{p:negative:residue:bound}
\end{equation}
for a suitable constant $c_5$. Moreover this bound is uniform in $V$.
Now, the rest of the argument follows \emph{mutatis mutandis} the same ideas as for the case $p>0$.
Namely, for any parameterisation $\Gamma_j(\sigma)$ of the segments $\Gamma_j$~\footnote[2]{We named it with the same letter since there is no problem of misunderstanding.} one can see that
$f(\Gamma_j(\sigma),\phi)=f(\Gamma_j(\sigma),0)+\mathcal{O}(\phi)$. So, taking into account the compactness of the $\sigma$-interval $[0,1]$, it follows that 
\begin{equation}
\int_{\Gamma_j} \frac{1}{f(y,p,\phi)} \, dy =
\int_0^1 \frac{1}{f(\Gamma_j(\sigma),\phi)} \, d\sigma =
\int_0^1 \left( \frac{1}{f(\Gamma_j(\sigma),0)} + \mathcal{O}(\phi) \right) \, d\sigma = c_{4+j} + \mathcal{O}(\phi),
\label{p:negative:Gammaj}
\end{equation}
for $j=2,3,4$. From expressions~\eqref{p:negative:time:1}, 
\eqref{p:negative:residue:bound}, and~\eqref{p:negative:Gammaj}
the following estimate for the flight time is derived:
\begin{equation}
T_{-}(\phi) = \tilde{C} +\mathcal{O}(\phi),   
\end{equation}
$\tilde{C}$ being a constant. In logarithmic scale this corresponds to say
$\log T_{-}(\phi) \simeq \mathrm{constant}$, a \textit{plateau}, different from the scaling law observed for $p>0$.


\begin{thebibliography}{99}
\bibliography{iopart-num}

\bibitem{Goldenfeld1992}
Goldenfeld N 
\newblock {\em Lectures On Phase Transitions And The Renormalization Group},
\newblock CRC Press 1992

\bibitem{Kuznetsov1998}
Kuznetsov Y 
\newblock {\em Elements of Applied Bifurcation Theory},
\newblock Second Edition. Springer 1998

\bibitem{Strogatz2000}
{Strogatz} SH (2000)
\newblock {\em Nonlinear Dynamics and Chaos with applications to
  Physics, Biology, Chemistry, and Engineering}.
\newblock Westview Press

\bibitem{Hohenberg1977}
Hohenber PC, Halperin BI 1977
\newblock The theory of dynamic critical phenomena.
\newblock {\em Rev. Mod. Phys.} {\bf 49} 435-479

\bibitem{Strogatz1989}
{Strogatz} SH, {Westervelt} RM 1989
\newblock {Predicted power laws for delayed switching of charge density waves}.
\newblock {\em Phys. Rev. B} {\bf 40} 10501-10508

\bibitem{Suzuki1982I}
Suzuki M, Kaneko K, Takesue S 1982
\newblock {Critical slowing down in stochastic processes I}.
\newblock {\em Progr. Theor. Phys.} {\bf 67} 1756-1775

\bibitem{Suzuki1982II}
Suzuki M, Takesue S, Sasagawa F 1982
\newblock {Critical slowing down in stochastic processes II}.
\newblock {\em Progr. Theor. Phys.} {\bf 68} 98-115

\bibitem{Sardanyes2019}
{Sardany\'es} J, {Pinero} J, {Sol\'e} R 2019 
\newblock {Habitat loss-induced tipping points in metapopulations with facilitation}.
\newblock {\em Pop Ecol} {\bf 61}(4) 436-449

\bibitem{Gimeno2018}
{Gimeno} J, {Jorba} \`A., {Sardany\'es}  2018
\newblock {On the effect of time lags on a saddle-node remnant in hyperbolic replicators}.
\newblock {J. Phys. A: Math. Theor.}, {\bf{51}} 385601

\bibitem{Fontich2008}
{Fontich} E, {Sardany\'es} J 2008
\newblock {General scaling law in the saddle-node bifurcation: a complex phase
  space study}.
\newblock {J. Phys. A: Math. Theor.} {\bf 41} 468-482.

\bibitem{Sardanyes2006}
{Sardany\'es} J, {Sol\'e} R 2006
\newblock {Bifurcations and phase transitions in spatially-extended two-member hypercycles.}
\newblock {\em J. Theor. Biol.} {\bf 234}(4) 468--482.

\bibitem{Canela2022}
J.~{Canela}, N.~{Fagella}, Ll.~{Alsed\`a}, and J.~{Sardany\'es} 2022
\newblock {Dynamical mechanism behind ghosts unveiled in a map complexification}.
\newblock {\it Chaos, Solitons \& Fractals} {\bf{156}} 111780. 

\bibitem{Trickey1998}
Trickey ST, Virgin LN 1998
\newblock Bottlenecking phenomenon near a saddle-node remnant in a {Duffing} oscillator.
\newblock {\em Phys. Lett. A}, {\bf 248} 185-190

\bibitem{Vidiella2018}
{Vidiella} B, {Sardany\'es} J, {Sol\'e} R 2018
\newblock {Exploiting delayed transitions to sustain semiarid ecosystems after catastrophic shifts.} 
\newblock {\em J. Roy. Soc. Interface} {\bf 15} 20180083.

\bibitem{Hill1910}
{Hill} AV 1910
\newblock { The combinations of haemoglobin with oxygen and with carbon monoxide.} 
\newblock {\em J. Physiol.} {\bf 40} iv-vii

\bibitem{SardanyesNJP2020}
Sardany\'es J, Raich C, Alarc\'on T 2020
\newblock{ Noise-induced stabilization of saddle-node ghosts.} 
\newblock {\em New J. Phys.} {\bf 22} 093064

\bibitem{Duarte2011}
Duarte J, Janu\'ario C, Martins N, Sardany\'es J 2012 
\newblock{Scaling law in saddle-node bifurcations for one-dimensional maps: a complex variable approach.}
\newblock {\em Nonlin. Dyn.} {\bf 67} 541-547

\bibitem{Kuehn2008}
{Kuehn} C 2008
\newblock {Scaling of saddle-node bifurcations: degeneracies and rapid quantitative changes.}
\newblock {\em J. Phys. A: Math Theor} {\bf 42}(4) 045101

\bibitem{Sardanyes2008}
{Sardany\'es} J 2008
\newblock {Error threshold ghosts in a simple hypercycle with error prone self-replication.}
\newblock {\em Chaos, Solitons \& Fractals} {\bf 35}(2) 313-319

\bibitem{Sardanyes2010}
{Sardany\'es} J, and {Fontich} E 2010
\newblock {On the metapopulation dynamics of autocatalysis: extinction transients related to ghosts.} 
\newblock {\em Int. J. Bifurcat. Chaos} {\bf 20}(4) 1-8

\bibitem{Sardanyes2007}
{Sardany\'es} J, {Sol\'e} R 2007
\newblock {The role of cooperation and parasites in non-linear replicator delayed extinctions.}
\newblock {\em Chaos, Solitons \& Fractals} {\bf 31} 1279-1296

\bibitem{assaf2006}
{Assaf} M, {Meerson} B 2006
\newblock {Spectral Theory of Metastability and Extinction in Birth-Death Systems.}
\newblock {\em Phys. Rev. E} {\bf 97} 200602

\bibitem{Assaf2010}
{Assaf} M, {Meerson} B 2010 
\newblock {Extinction of metastable sto-chastic populations.}
\newblock {\em Phys. Rev. E} {\bf 81} 021116

\bibitem{escudero2009}
{Escudero} C, {Kamenev} A 2009 
\newblock {Switching rates of multistep reactions.}
\newblock {\em Phys. Rev. E} {\bf 79} 041149

\bibitem{gottesman2012}
{Gottesman} O, {Meerson} B 2012 
\newblock {Multiple extinction routes in stochastic population models.}
\newblock {\em Phys. Rev. E} {\bf 85} 021140

\bibitem{bressloff2014a}
Bressloff PC
\newblock {\em Stochastic processes in cell biology.},
\newblock Springer-Verlag, Berlin, Germany 2014

\bibitem{hinch2005}
{Hinch} R, {Chapman} SJ 2005 
\newblock {Exponentially slow transitions on a Markov chain: the
frequency of calcium sparks.}
\newblock {\em Eur. J. Appl. Math.} {\bf 16} 427-446

\bibitem{elgart2004}
{Elgart} V, {Kamenev} A 2004
\newblock {Rare event statistics in reaction-diffusion systems.}
\newblock {\em Phys. Rev. E} {\bf 70} 041106


\bibitem{Gillespie1977}
{Gillespie} DT 1977
\newblock {Exact stochastic simulation of coupled chemical reactions}.
\newblock {\em J. Phys. Chem.}, {\bf{81}}(25) 2340-2436

\bibitem{Gillespie2000}
{Gillespie} DT 2000
\newblock {The chemical Langevin equation.}
\newblock {\em J. Chem. Phys.}, {\bf{113}}: 297-306

\bibitem{Kubo1973}
Kubo R, Matsuo K,  Kitahara K 1973 
\newblock {Fluctuation and Relaxation of Macrovariables.}
\newblock {\em J. Stat. Phys.} {\bf{9}}(1) 51-96 

\end{thebibliography}
\end{document}